\documentclass{amsart}

\usepackage{amssymb}
\usepackage{amsmath}
\usepackage{array}

\usepackage{wasysym}

\usepackage{listings}

\newcommand{\Pf}{{\em Proof}. }
\newcommand{\EPf}{\hfill$\square$}

\newcommand{\Lg}{\mbox{$\mathfrak g$}}

\newcommand{\Lh}{\mbox{$\mathfrak h$}}
\newcommand{\Lk}{\mbox{$\mathfrak k$}}

\newcommand{\Ln}{\mbox{$\mathfrak n$}}

\newcommand{\Lt}{\mbox{$\mathfrak t$}}

\newcommand{\Lz}{\mbox{$\mathfrak z$}}

\newcommand\liegr{\sf}

\newcommand{\SO}[1]{\mbox{${\liegr SO}(#1)$}}

\newcommand\fieldsetc{\mathbb}

\newcommand{\R}{\fieldsetc{R}}
\newcommand{\C}{\fieldsetc{C}}

\newtheorem*{thm}{Theorem}

\newtheorem{prop}{Proposition}
\newtheorem{lem}{Lemma}

\newcommand{\sff}{\mathrm{II}}

\title{Focal radii of orbits}

\author{Claudio Gorodski and Artur B. Saturnino}

\address{Instituto de Matem\'atica e Estat\'\i stica, Universidade de
S\~ao Paulo, Rua do Mat\~ao, 1010, S\~ao Paulo, SP 05508-090, Brazil}

\email{gorodski@ime.usp.br}

\address{Department of Mathematics, University of Pennsylvania,
  209 South 33rd Street,
Philadelphia, PA 19104-6395, USA}

\email{bsatur@sas.upenn.edu}

\thanks{The first author was partially supported by the
  CNPq grant 303038/2013-6 and the FAPESP project 2011/21362-2.}

\thanks{The second author was supported by the CNPq scholarship
  154010/2015-4.}


\begin{document}

\begin{abstract}
We show that every effective action of a compact Lie group $K$
on a unit sphere $S^n$ admits an explicit orbit whose principal 
curvatures are bounded from above by~$4\sqrt{14}$. 
\end{abstract}

\maketitle

Let $K$ be a compact Lie group acting by isometries on the 
unit sphere $S^n$. The orbit space $X=S^n/K$ is an Alexandrov
space of curvature bound below by $1$ and diameter bounded 
above by $\pi$. In this context, the \emph{range problem}
of Grove and Markovsen~\cite{GM} can take the following 
interpretation:
\begin{quote}
\emph{How small can the diameter of $X$ be?}
\end{quote}
Assume the action of $K$ on $S^n$ is nontransitive and $n\geq2$
to dismiss trivial cases. The best result that has been achieved
so far is the existence of a non-explicit, dimension-dependent
positive lower bound for the diameter of~$X$ that however goes to zero 
as $n\to\infty$~\cite{Gr}. The analysis of several special
cases~\cite{McG,Gr,Searlediameter} 
and a related result on the curvature of $X$~\cite{GL4} 
indeed suggest that a dimension-independent positive lower
bound should exist. 

For an arbitrary metric space $X$, recall that
the \emph{diameter} is
\[ \mathrm{diam}(X)=
\sup\{\,d(x,y)\;|\;x,\ y\in X\,\}\]
and the \emph{radius} at $x\in X$ 
is $r_x=\inf\{\,r>0\;|\;X\subset B(x,r)\,\}$. It is immediate from
the triangle inequality that 
\begin{equation*}
 r_x \leq \mathrm{diam}(X) \leq 2r_x 
\end{equation*}
for all $x\in X$. It follows from these inequalities that bounding the 
diameter of $X$ is equivalent to bounding the radius of $X$ at an
arbitrary, fixed point $x\in X$. 

Let $M$ be a $K$-orbit in $S^n$ represented by a point $x\in X=S^n/K$.
Recall that a focal point of $M$ is 
simply a critical value of the 
normal exponential map of $M$ in $S^n$. Define
$\overline{f}_M$, $\underline{f}_M$ to be the supremum, resp.\ the infimum
of the focal distances to $M$ along normal geodesics.
While it is reasonable to call $\underline{f}_M$ the \emph{focal radius} of $M$,
the other number is also of interest to us since
$\overline{f}_M\geq r_x\geq\frac12\mathrm{diam}(X)$, and we see 
that the existence of a uniform lower bound for the maximal focal distance
of an arbitrary $K$-orbit is a necessary condition for the existence of a 
uniform lower bound for the diameter of $X$. Our main result 
confirms this necessary condition for a specific $K$-orbit but in
the stronger sense of bounding the focal radius $\underline{f}_M$.

In order to state our result, 
we assign a number $C_{\mathbf T}$
to each Cartan type $\mathbf T$
of simple Lie groups according to the following table:
\[\begin{array}{|c|c|c|c|c|c|c|c|c|}
\hline
A_n & B_n & C_n & D_n & G_2 & F_4 & E_6 & E_7 & E_8 \\
\hline
2\sqrt3 & 2\sqrt6 & 2\sqrt6 & 2\sqrt6 & \sqrt6 & 2\sqrt7 & 2\sqrt{10}
& 8 & 4\sqrt 7\\
\hline
\end{array}
\]
\begin{center}
  \textsc{Table~1}
\end{center}
We only consider the case in which
the associated representation of $K$ on $\R^{n+1}$ is irreducible
for the sake of simplicity and owing to the fact that 
otherwise $\mathrm{diam}(X)\geq\pi/2$~\cite{GL}. 

\begin{thm}
Let $X$ be the orbit space of an effective action of a compact Lie group $K$
on the unit sphere $S^n$. Assume that the associated 
orthogonal representation $\rho$ of $K$ on $\R^{n+1}$ is irreducible.  
Let $C$ be the maximum of the numbers $C_{\mathbf T}$, 
where $\mathbf T$ runs through the Cartan types  
of the simple factors of $K$. Then the focal radius of a certain
$K$-orbit in $S^n$ is bounded below by $\mathrm{arccot}(C)$
or $\mathrm{arccot}(C\sqrt2)$, according to whether $\rho$ is 
admits an invariant complex structure or not. 
\end{thm}

Understanding the diameter of quotients of unit spheres is an interesting
problem on its own, but it also has a bearing on the global structure 
and classification of compact positively curved manifolds.
In fact the study of manifolds of positive 
curvature usually starts with
those with large groups of symmetries, and one is naturally 
led to consider the corresponding orbit spaces. The orbit space of a 
compact positively curved Riemannian manifold under the action of a
compact Lie group is an Alexandrov space of positive curvature. 
The local geometry of such an Alexandrov space is given by 
its tangent cones which, in turn, are described in terms of 
orbit spaces of unit spheres under actions of compact Lie groups  
(see e.g.~\cite[Theorem~1.4]{GroveSearle} 
and~\cite[Theorem~2.5]{wilking-survey} for 
concrete examples).

The orbit space $X=S^n/K$ can also be viewed as the bi-quotient 
$\SO{n+1}/\!/(K\times\SO n)$, where $K$ acts on $\SO{n+1}$ from 
the left and $\SO n$ acts from the right. Let $G$ be a compact connected
Lie group equipped with a bi-invariant Riemannian metric. Then 
$G\times G$ acts by isometries on $G$ via left and right translations.
For a closed subgroup $H$ of $G\times G$, one can consider the 
bi-quotient $G/\!/H$ and ask
\begin{quote}
\emph{How small can the diameter of $G/\!/H$ be?}
\end{quote}

The authors would like to thank Francisco Gozzi, Alexander Lytchak 
and Luiz San Martin for informative discussions.
This work is based on the Masters Dissertation of the
second named author.

\section{Preliminaries}

\subsection{Geometric remarks}

Since $X=S^n/K$ is compact, 
the diameter and the radii are realized by minimizing geodesics. 
The horizontal lift of a minimizing geodesic in $X$ is a 
minimizing geodesic in $S^n$ between two $K$-orbits,
orthogonal to them, and of the same length. Fix an orbit
$M=Kp$. A unit speed
geodesic $\gamma:[0,+\infty)\to S^n$, starting orthogonally 
to $M$ at $p$, can only cease to be minimizing at $\gamma(t_0)$ for some $t_0>0$
for one of the following 
two non-mutually exclusive reasons: (i) $\gamma(t_0)$ is the first focal point
of $M$ along $\gamma$; (ii) there is another minimizing geodesic 
from $M$ to $\gamma(t_0)$ of the same length as $\gamma$. 
The condition~(i) is a 
local one and, as we shall see, within control for a certain choice of $p$. 
Indeed, since the ambient space $S^n$ 
has constant curvature $1$, $\gamma(t_0)$ is the first focal point
of $M$ along $\gamma$ if and only if $\cot t_0$ is the largest 
eigenvalue of the Weingarten operator $A_\xi$ of $M$, where $\xi=\gamma'(0)$. 
Denote the second fundamental form of $M$ at~$p$ by
$\mathrm{II}_p:T_pM\times T_pM\to\nu_pM$. If $\cot t_0$ is 
the largest eigenvalue of $A_\xi$, then we have the following
inequalities for the supremum norm:
\[ ||\mathrm{II}_p||_\infty\geq||A_\xi||_\infty=\cot t_0\geq\cot\underline{f}_M. \]
We deduce that the focal radius of $M$ satisfies
\begin{equation}\label{focal-radius}
 \underline{f}_M\geq\mathrm{arccot}(||\mathrm{II}_p||_\infty). 
\end{equation}

\subsection{Algebraic remarks}\label{algebraic}
Denote by $\Lk$ the Lie algebra of $K$. Choose a maximal torus $T$
of $K$ with associated Lie algebra $\Lt$. Then the complexification
$\Lh:=\Lt^{\mathbb C}$ is a Cartan subalgebra of $\Lg:=\Lk^{\mathbb C}$
and we have the corresponding root decomposition 
$\Lg=\Lh+\sum_{\alpha\in\Delta}\Lg_\alpha$. Choose an ordering of the roots
and denote by $\Delta^+$ the corresponding system of positive roots.
We will denote the Cartan-Killing form of $\Lg$ (and its transfer
to the dual~$\Lg^*$) by $\langle\cdot,\cdot\rangle$.

It is known that we can choose root vectors $e_\alpha\in\Lg_\alpha$
and $f_\alpha\in\Lg_{-\alpha}$ for each $\alpha\in\Delta^+$ so that
\begin{equation}\label{1}
 [e_\alpha,f_\alpha]=h_\alpha 
\end{equation}
where
\[ \langle h_\alpha,h\rangle = \alpha^\vee(h)=\frac{2\alpha(h)}{||\alpha||^2} \]
for all~$h\in\Lh$, and
\[ x_\alpha:=e_\alpha-f_\alpha,\ y_\alpha:=i(e_\alpha+f_\alpha),\ ih_\alpha \]
span the semisimple part of $\Lk$.

\subsubsection{Complex representations}

Assume $\rho$ admits an invariant complex structure
and view it as a representation 
of $K$ on a complex vector space $V$. It is a standard fact that
the complex structure on $V$ is an orthogonal transformation,
so there is a canonical extension of the $K$-invariant (real) inner product
to a $K$-invariant
Hermitian product on $V$ which we also shall denote 
by $\langle\cdot,\cdot\rangle$, without risking ambiguity. 
Finally, $V$ can also be regarded as a representation of $\Lg$. 

\begin{lem}\label{norm}
  Let $\alpha\in\Delta^+$.
  \begin{enumerate}
  \item 
    The adjoint of the operator $e_\alpha$ on $V$ is $f_\alpha$,
    namely,
    \[ \langle e_\alpha v_1,v_2\rangle=\langle v_1,f_\alpha v_2\rangle\]
    for all $v_1$, $v_2\in V$.
  \item If $v_\mu$ is a weight vector of $V$ of weight $\mu$, then
    \[ ||f_\alpha v_\mu||^2=||e_\alpha v_\mu||^2+ \langle \mu,\alpha^\vee\rangle||v_\mu||^2. \]
    \end{enumerate}
    \end{lem}

\Pf Part~(a) is well known to follow from the fact that the elements
of $\Lk$ act on $V$ by skew-adjoint endomorphisms. For part~(b), we 
compute 
\begin{eqnarray*}
\langle f_\alpha v_\mu,f_\alpha v_\mu\rangle&=&\langle e_\alpha f_\alpha v_\mu,v_\mu\rangle \\
&=& \langle f_\alpha e_\alpha v_\mu,v_\mu\rangle +\langle h_\alpha v_\mu,v_\mu\rangle\\
&=& \langle e_\alpha v_\mu,e_\alpha v_\mu\rangle+\langle\mu,\alpha^\vee\rangle
\langle v_\mu,v_\mu\rangle,
\end{eqnarray*}
where we have used part~(a) and equation~(\ref{1}). \EPf

\subsection{Isotropy representations of symmetric spaces}

It is convenient to have the lemma below for later use. 

\begin{lem}\label{isoparametric}
  If $\rho$ is the isotropy representation of a symmetric space
  of rank greater than one,
  then we can choose $p\in S^n$
  such that $||\sff_p||\leq2$.
\end{lem}

\Pf In this case the $K$-orbits comprise an isoparametric foliation
of $\R^{n+1}$ of codimension at least two and the principal curvatures
of the leaves are known explicitly. Namely, let $\Sigma$ be a
fixed normal space to a principal $K$-orbit. Then $\Sigma$ meets
every $K$-orbit.
There is a reduced root system $\Delta$ in the dual space
$\Sigma^*$ such that the principal curvatures of $Kp$ for $p\in\Sigma$ 
are given by $-\alpha(\xi)/\alpha(p)$ for all $\xi\in\Sigma\subset\nu_p(Kp)$
and $\alpha\in\Delta$ satisfying $\alpha(p)\neq0$~\cite[Example~2.7.1]{bco}.

Let $\tilde\alpha$ be the highest root of $\Delta$ with respect to some
ordering of the roots and take $p= \frac{||\tilde\alpha||}{2}h_{\tilde\alpha}$, where
$h_{\tilde\alpha}$ is as in Subsection~\ref{algebraic}. Then, considering a unit
vector $\xi$ and $\alpha$ non-orthogonal to $\tilde\alpha$, we have
\[ \left|\frac{\alpha(\xi)}{\alpha(p)}\right|\leq
\frac{||\alpha||\, ||\tilde\alpha||}{|\langle\alpha,\tilde\alpha\rangle|}
=\frac1{|\cos\theta|}\leq2, \]
where $\theta$ is the angle between $\alpha$ and $\tilde\alpha$,
by the cristallographic property of root systems. \EPf

\section{Setting for case of complex representations}

Our standing assumption in this section and 
in Sections~\ref{reduction} and~\ref{classical}
is that~$\rho$ leaves a complex
structure invariant. 
Then it can be considered as complex representation~$V$
and we can apply Cartan's highest weight theory. 

Denote by $\lambda$ the highest weight of $V$ and choose a unit highest
weight vector $v_\lambda$. Then $p=v_\lambda$ is a point in $S^n$.
Consider the orbit $M=Kp$. 
Our goal is to estimate the supremum norm of the second fundamental 
form of the orbit $Kp$ (at $p$):
\[ ||\Pi_p||_\infty =\sup\{\,||\Pi(u,v)||\;|\;||u||=||v||=1\,\}. \]

The center $\Lz$ of $\Lk$ acts as a scalar 
on $V$, so the tangent space
\begin{eqnarray*} 
T_pM &=& \Lk p\\
    &=& \Lz p+\sum_{\alpha\in\Delta^+}\R ih_\alpha p+ \R x_\alpha p+ \R y_\alpha p \\
&=& i\R p + \sum_{\alpha\in\Delta^+} \C f_\alpha p \\
&=& i\R p + \Ln_- p
\end{eqnarray*}
where $\Ln_-=\sum_{\alpha\in\Delta^+}\Lg_{-\alpha}$. 
Then $\C p +\Ln_-p$ is a complex subspace of $V$
and its (Hermitian) orthogonal complement in $V$ is the normal space
to $M$ at~$p$ in $S^n$:
\[ \nu_p M = (\C p +\Ln_- p)^\perp. \]

The second fundamental form
\[ \sff_p:T_pM\times T_pM \to \nu_p M \]
is given by 
\[ \sff_p(xp,yp)=(xyp)^\nu \]
for $x$, $y\in\Lk$, where $(\cdot)^\nu$ denotes the 
component in $\nu_pM$. 

It is obvious that $\sff_p(i\R p,T_pM)=0$, so it suffices to
consider the restriction of $\sff_p$ to $\Ln_-\times\Ln_-$. 

\begin{lem}\label{complex}
The restriction 
\[ \sff_p:\Ln_-p\times\Ln_-p\to\nu_pM \]
is $\C$-bilinear. It follows that 
\[ \sff_p(fp,f'p)=(ff'p)^\nu \]
for all $f$, $f'\in\Ln_-$. 
\end{lem}

\Pf For the first assertion,
by symmetry it suffices to show that $\sff_p$ is $\C$-linear
in the second argument. Let $x$, $y\in\Lk$ be such that $xp$,
$yp\in\Ln_-p$. Then $i(yp)\in\Ln_-p$ so there exists $\tilde y\in\Lk$
such that $i(yp)=\tilde yp$. Now
\begin{eqnarray*}
\sff_p(xp,i(yp))&=&\sff_p(xp,\tilde yp)\\
  &=& (x\tilde yp)^\nu \\
  &=& (xi(yp))^\nu \\
  &=& (ix(yp))^\nu \\
  &=& i(xyp)^\nu \\
  &=& i\sff_p(xp,yp),
\end{eqnarray*}
as wished. The second assertion follows from the first one and 
the fact that the orthogonal projection to $\nu_pM$ is $\C$-linear. \EPf

\subsection{Reduction to case $K$ is simple}\label{sec:simple}

By passing to an almost effective action, we may assume that 
$K$ splits as the direct sum of a torus (its center) and 
its simple factors $K_1,\ldots, K_m$. Since $V$ is complex irreducible,
the center is at most one-dimensional and $V$ can be written as a complex tensor product $V_1\otimes\cdots\otimes V_m$ 
where $V_i$ is a complex irreducible representation of $K_i$.
We will show that it suffices to consider the representations of 
$K_i$ on $V_i$. 

Recall that $p=v_\lambda$ is a unit highest weight vector of $V$.
Note that $K$ having a center or not is irrelevant for the $K$-orbit
through~$p$ and for the irreducibility of $V$. 
We can write 
$p=p_1\otimes\cdots\otimes p_m$ where $p_i$ is a 
unit highest weight vector of $V_i$. We will use that the 
Hermitian product on $V$ is given in terms of the 
Hermitian products on the $V_i$ by
\[ \langle u_1\otimes\cdots\otimes u_m,
u'_1\otimes\cdots\otimes u'_m\rangle = \langle u_1,u'_1\rangle
\cdots\langle u_m,u'_m\rangle, \]
and that
\[ ||\sff_p||_\infty=\sup\{||\sff(fp,fp)||;f\in\Ln_-,\ ||fp||=1\} \]
since $\sff_p$ is symmetric.

\begin{prop}\label{simple}
Let $\sff$ denote the second fundamental form of $Kp$ at~$p$ in the unit 
sphere of $V$ and 
let $\sff_i$ denote the second fundamental form of $K_ip$ at~$p_i$ in the 
unit sphere of $V_i$. If $C\geq\max\{||\sff_1||_\infty,\ldots,||\sff_m||_\infty,
\sqrt2\}$ then $C\geq||\sff||_\infty$. 
\end{prop}

\Pf Let $\Lg_i$ be the complexification of the Lie algebra $\Lk_i$  
of $K_i$, $\Ln_{i-}$ the corresponding negative nilpotent subalgebra
and $\Ln_-=\sum_{i=1}^m\Ln_{i-}$. Let $f=\sum_{i=1}^mf_i$ where $f_i\in\Ln_{i-}$. 

We first remark that if $i<j$ then 
\begin{eqnarray*}
||\sff(f_ip,f_jp)||^2&=&||(f_if_jp)^\nu||^2\\
&\leq&||f_if_jp||^2\\
&=&||p_1\otimes \cdots\otimes f_ip_i\otimes\cdots\otimes f_jp_j\otimes
\cdots\otimes  p_m||^2\\
&=& ||f_ip_i||^2||f_jp_j||^2\\
&=& ||f_ip||^2||f_jp||^2.
\end{eqnarray*}

Note also that
\begin{eqnarray*}
||\sff(f_ip,f_ip)||^2&=&||(p_1\otimes\cdots\otimes f_i^2p_i\otimes\cdots\otimes p_m)^\nu||^2 \\
&=&||p_1\otimes\cdots\otimes (f_i^2p_i)^\nu\otimes\cdots\otimes p_m||^2 \\
&=&||(f_i^2p_i)^\nu||^2\\
&=&||\sff_i(f_ip_i,f_ip_i)||^2\\
&\leq& C^2||f_ip_i||^4\\
&=& C^2||f_ip||^4.
\end{eqnarray*}

Finally, owing to the fact that
weight spaces associated to different weights are 
orthogonal, and using the remarks above, we can write
\begin{eqnarray*}
||\sff(fp,fp)||^2& = & ||\sum_{i,j=1}^m\sff(f_ip,f_jp)||^2 \nonumber \\
&=& \sum_{i=1}^m||\sff(f_ip,f_ip)||^2+ 4\sum_{i<j}||\sff(f_i,f_j)||^2 \nonumber \\
&\leq& C^2\sum_{i=1}^m||f_ip||^4+ 4\sum_{i<j}||f_ip||^2||f_jp||^2\nonumber \\
&\leq&C^2\left(\sum_{i=1}^m||f_ip||^4+ 2\sum_{i<j}||f_ip||^2||f_jp||^2\right) \nonumber\\
&=&C^2||fp||^4,
\end{eqnarray*}
as desired. \EPf

\subsection{Estimate on basis vectors}

Consider the complex basis of $\Ln_-p$ given by the 
$f_\alpha p$ for $\alpha\in\Delta^+$. 

\begin{lem}\label{mab}
For all $\alpha$, $\beta\in\Delta^+$ we have
\[ ||\sff_p(f_\alpha p,f_\beta p)||\leq m_{\alpha,\beta}\,||f_\alpha p||\,||f_\beta p||
\]
where 
\[ m_{\alpha,\beta}=\left\{\begin{array}{rl}\sqrt2&\mbox{if $\alpha=\beta$ or
$\langle \alpha,\beta\rangle<0$,}\\
1&\mbox{otherwise.}\end{array}\right. \]
\end{lem}

\Pf Since $p$ is a unit highest weight vector, it follows from 
Lemma~\ref{norm}(b) that 
\[ ||f_\alpha p||^2= \langle\lambda,\alpha^\vee\rangle\quad\mbox{and}\quad
||f_\beta p||^2= \langle\lambda,\beta^\vee\rangle. \]
Without loss of generality, we may assume these numbers are not zero. 

In view of Lemma~\ref{complex}
and using Lemma~\ref{norm} with $\mu=\lambda-\beta$, we get
\begin{eqnarray*}
||\sff_p(f_\alpha p,f_\beta p)||^2 &=& ||(f_\alpha f_\beta p)^\nu||^2\\
&\leq& ||f_\alpha f_\beta p||^2\\
&=& || e_\alpha f_\beta p ||^2 + \langle \lambda-\beta,\alpha^\vee\rangle||f_\beta p||^2\\
&=&||[e_\alpha,f_\beta]p||^2+||f_\alpha p||^2||f_\beta p||^2-\langle\beta,\alpha^\vee
\rangle||f_\beta p||^2.
\end{eqnarray*}
We next consider different cases.

1.~$\alpha=\beta$. Here
\begin{eqnarray*}
||\sff_p(f_\alpha p,f_\alpha p)||^2 &=&
||h_\alpha p||^2+||f_\alpha p||^4-2||f_\alpha p||^2\\
& \leq & \langle\lambda,\alpha^\vee\rangle^2+||f_\alpha p||^4\\
&=& 2||f_\alpha p ||^2.
\end{eqnarray*}

2.~$\alpha\neq\beta$ and $\langle\alpha,\beta\rangle\geq0$. 
Since $\sff_p$ is symmetric, changing the roles of $\alpha$ and $\beta$
we may assume that $\beta-\alpha\not\in\Delta^+$. Then 
$[e_\alpha,f_\beta]p=0$ and we get 
\[ ||\sff_p(f_\alpha p,f_\beta p)||^2
=||f_\alpha p||^2||f_\beta p||^2-\langle\beta,\alpha^\vee
\rangle||f_\beta p||^2\leq
 ||f_\alpha p||^2||f_\beta p||^2.\]

3.~$||\alpha||=||\beta||$ and $\langle\alpha,\beta\rangle<0$. 
As in case 2, we may assume $\beta-\alpha\not\in\Delta^+$ and 
thus $[e_\alpha,f_\beta]p=0$. By~\cite[ch.~VI, \S1, no.3]{Bourbaki},
we have $\langle\beta,\alpha^\vee\rangle=-1$. Owing to
$||f_\alpha p||^2=\langle\lambda,\alpha^\vee\rangle\geq1$, we deduce
 \[ ||\sff_p(f_\alpha p,f_\beta p)||^2
=||f_\alpha p||^2||f_\beta p||^2-\langle\beta,\alpha^\vee
\rangle||f_\beta p||^2\leq
2||f_\alpha p||^2||f_\beta p||^2.\]

4.~$||\alpha||\neq||\beta||$ and $\langle\alpha,\beta\rangle<0$. 
Changing the roles of $\alpha$ and $\beta$, we may assume that 
$\alpha$ is long and $\beta$ is short. By~\cite[ch.~VI, \S1, no.3]{Bourbaki},
we get $\langle\beta,\alpha^\vee\rangle=-1$. 
Since $\langle\alpha,\beta\rangle<0$, we have $||\beta-\alpha||>||\alpha||$
and thus $\beta-\alpha$ cannot be a root, for otherwise $\alpha$ and $\beta$
would span an irreducible root subsystem of $\Delta$ with three different
root lengths. We deduce that  $[e_\alpha,f_\beta]p=0$ and hence
 \[ ||\sff_p(f_\alpha p,f_\beta p)||^2
=||f_\alpha p||^2||f_\beta p||^2-\langle\beta,\alpha^\vee
\rangle||f_\beta p||^2\leq
2||f_\alpha p||^2||f_\beta p||^2\]
as in case~3. \EPf

\section{Reduction to a problem about root systems}\label{reduction}

Let $f\in\Ln_-$ be such that $||fp||=1$ and
write $f=\sum_{\alpha\in\Delta^+}z_\alpha f_\alpha$ for 
$z_\alpha\in\C$. To avoid dividing by zero, we put
\[ r_\alpha = |z_\alpha|\,||f_\alpha p|| \]
for all $\alpha\in\Delta^+$; note that 
\[ \sum_{\alpha\in\Delta^+} r_\alpha^2 = \sum_{\alpha\in\Delta^+}|z_\alpha|^2||f_\alpha p||^2
=||fp||^2=1. \]

Using the fact that the orthogonal projection onto $\nu_pM$ preserves weight spaces and applying Lemma~\ref{mab}, we can write
\begin{eqnarray} \nonumber
||\sff_p(fp,fp)||^2&=&\left\lVert\sum_{\alpha,\beta\in\Delta^+}z_\alpha z_\beta\sff(f_\alpha p,f_\beta p)\right\rVert^2\\ \nonumber
&=&\sum_{\gamma\in2\Delta^+}\left\lVert\sum_{\alpha+\beta=\gamma}z_\alpha z_\beta\sff(f_\alpha p,f_\beta p)\right\Vert^2\\ \nonumber
&\leq&\sum_{\gamma\in2\Delta^+}\left(\sum_{\alpha+\beta=\gamma}||z_\alpha z_\beta\sff(f_\alpha p,f_\beta p)||\right)^2\\ \label{partial-estimate}
&\leq&\sum_{\gamma\in2\Delta^+}\left(\sum_{\alpha+\beta=\gamma}m_{\alpha,\beta}r_\alpha r_\beta
\right)^2.
\end{eqnarray}
At this juncture, we can use the Cauchy-Schwarz inequality to deduce
\begin{eqnarray} \nonumber
  ||\sff_p(fp,fp)||^2&\leq&\sum_{\gamma\in2\Delta^+}\left(\sum_{\alpha+\beta=\gamma}m_{\alpha,\beta}^2\right)
  \left(\sum_{\alpha+\beta=\gamma}r_\alpha^2r_\beta^2\right)\\ \nonumber
  &\leq& C_\Delta\sum_{\gamma\in2\Delta^+}
\left(\sum_{\alpha+\beta=\gamma}r_\alpha^2r_\beta^2\right) \\ \nonumber
&=&C_\Delta\left(\sum_{\alpha\in\Delta^+} r_\alpha^2\right)^2\\ \label{weak-estimate}
&=&C_\Delta, 
\end{eqnarray}
where we have defined 
\[ C_\Delta= \max_{\gamma\in2\Delta^+}\sum_{\alpha+\beta=\gamma}m_{\alpha,\beta}^2.\]

\subsection{The case $K$ is a simple Lie group
  of exceptional type}\label{exceptional}

The estimate~(\ref{weak-estimate}) is relatively easy to obtain, but not enough
for our purposes since $C_\Delta\to\infty$ as the rank of $K$ 
goes to infinity. Nonetheless, due to Proposition~\ref{simple},
we may assume that $K$ is a 
simple Lie group and then we can use this estimate in case $K$ is of
exceptional type. Using a simple algorithm to compute~$C_\Delta$
(see App.), we deduce that $C_\Delta=\sqrt6$, $2\sqrt7$, $2\sqrt{10}$,
$8$ or $4\sqrt7$ according to whether $K$ is of type~$G_2$, $F_4$,
$E_6$, $E_7$ or $E_8$.

\section{The case $K$ is a classical Lie group}\label{classical}

We keep the notation from Section~\ref{reduction}.
Starting from~(\ref{partial-estimate}), we obtain
a more refined estimate on $||\sff_p||_\infty$ for each classical
family of compact simple Lie groups as follows. First note that
\[ \sum_{\gamma\in2\Delta^+}\sum_{\alpha+\beta=\gamma}r_\alpha^2r_\beta^2
=\sum_{\alpha,\beta\in\Delta^+}r_\alpha^2r_\beta^2=||fp||^4=1, \]
so we can rewrite~(\ref{partial-estimate}) as
\begin{equation}\label{II}
  ||\sff_p(fp,fp)||^2\leq 2+ \sum_{\gamma_\in2\Delta^+}S_\gamma
  \end{equation}
where
\[ S_\gamma = \left(\sum_{\alpha+\beta=\gamma}m_{\alpha,\beta}r_\alpha r_\beta\right) ^2-
2\sum_{\alpha+\beta=\gamma}r_\alpha^2r_\beta^2 \]
for all $\gamma\in2\Delta^+$. The next lemma explains
why we have singled out the $S_\gamma$.

\begin{lem}\label{sgamma}
  We have $S_\gamma>0$ for some $\gamma\in2\Delta^+$ only if
  $\gamma$ can be written as a sum of two
  positive roots in more than two ways
  (counting permutations) or $\gamma$ is the sum of two positive roots
  forming an obtuse angle.
\end{lem}

\Pf Suppose there is a unique way of writing $\gamma=\alpha+\beta$
with $\alpha$, $\beta\in\Delta^+$. Then $\alpha=\beta$ and $m_{\alpha,\alpha}=\sqrt2$,
so $S_\gamma=(\sqrt2 r_\alpha^2)^2-2r_\alpha^4=0$. Suppose next that
there are exactly two ways of decomposing $\gamma$, namely,
$\gamma=\alpha+\beta=\beta+\alpha$, and $\langle\alpha,\beta\rangle\geq0$.
Then $m_{\alpha,\beta}=m_{\beta,\alpha}=1$ and
\[ S_\gamma=(r_\alpha r_\beta+r_\beta r_\alpha)^2-2(r_\alpha^2 r_\beta^2+r_\beta^2 r_\alpha^2)=
0, \]
as desired. \EPf

\medskip

We shall estimate $\sum_{\gamma\in2\Delta^+}S_\gamma$ for each classical
family of simple Lie groups. Due to Lemma~\ref{sgamma}, 
$\sum_{\gamma\in2\Delta^+}S_\gamma=\sum_{\gamma\in\Phi}S_\gamma$,
where $\Phi$ is the subset
of $2\Delta^+$ consisting of elements $\gamma$
that satisfy the following condition:
\emph{$\gamma$ can be written as a sum of two positive roots in more
than two ways (counting permutations) or $\gamma$ is the sum of two positive roots
forming an obtuse angle.}

\subsection{The $A_n$ family}
Here $\Phi=\Phi_1\cup\Phi_2$ (disjoint union) where
\[ \Phi_1=\{\theta_i+\theta_j-\theta_k-\theta_\ell|i<j<k<\ell\} \]
and
\[ \Phi_2=\{\theta_i-\theta_j|i+1<j\}. \]
We will prove that
\[ \sum_{\gamma\in\Phi_1}S_\gamma\leq2\quad\mbox{and}\quad
\sum_{\gamma\in\Phi_2}S_\gamma\leq8. \]
It will thus follow from~(\ref{II}) that $||\sff_p||_\infty\leq2\sqrt3$. 

The trick is to consider the $(n+1)\times(n+1)$ matrix $A=(a_{ij})$
given by
\[ a_{ij}=\left\{\begin{array}{cl}r_{\theta_i-\theta_j}&\mbox{if $i<j$,}\\
0 &\mbox{otherwise.}
\end{array}\right. \]
Note that the Euclidean norm $||A||^2=\sum_{i<j}r_{\theta_i-\theta_j}^2=1$. 
We will also use the Euclidean inner product $\langle\cdot,\cdot\rangle$
on the space of  $(n+1)\times(n+1)$ real matrices.

If $\gamma=\alpha+\beta\in\Phi_1$ with $\alpha$, $\beta\in\Delta^+$ then, up to
permuting $\alpha$ and $\beta$, we must have 
$\alpha=\theta_i-\theta_k$ and $\beta=\theta_j-\theta_\ell$ or
$\alpha=\theta_i-\theta_\ell$ and $\beta=\theta_j-\theta_k$,
where $i<j<k<\ell$. In both cases $\langle\alpha,\beta\rangle=0$ so
$m_{\alpha,\beta}=1$ and therefore
\[ S_{\theta_i+\theta_j-\theta_k-\theta_\ell}=(2a_{ik}a_{j\ell}+2a_{i\ell}a_{jk})^2
-2(2a_{ik}^2a_{j\ell}^2+2a_{i\ell}^2a_{jk}^2)=8a_{ik}a_{j\ell}a_{i\ell}a_{jk}. \]
Since the expression above is symmetric in~$i$, $j$ and in $k$, $\ell$,
we can write
\begin{eqnarray*}
  \sum_{\gamma\in\Phi_1}S_\gamma&=&8\sum_{i<j<k<\ell}a_{ik}a_{j\ell}a_{i\ell}a_{jk}\\
  &=&2\sum_{i,j<k,\ell}a_{ik}a_{j\ell}a_{i\ell}a_{jk}\\
  &\leq&2\sum_{i,j,k,\ell}a_{ik}a_{j\ell}a_{i\ell}a_{jk}\\
  &=&2\langle AA^t,AA^t\rangle\\
  &\leq&2||AA^t||^2\\
  &\leq&2||A||^4\\
  &=&2.
\end{eqnarray*}

If $\gamma=\alpha+\beta\in\Phi_2$ with $\alpha$, $\beta\in\Delta^+$ then,
up to permuting $\alpha$
and $\beta$, we must have $\alpha=\theta_i-\theta_k$ and
$\beta=\theta_k-\theta_j$ with $i<k<j$. In this situation
$m_{\alpha,\beta}=\sqrt2$ so
\[ S_{\theta_i-\theta_j}=\left(\sum_{k=i+1}^{j-1}2\sqrt2a_{ik}a_{kj}\right)^2
-2\sum_{k=i+1}^{j-1}2a_{ik}^2a_{kj}^2\le 8\left(\sum_k a_{ik}a_{kj}\right)^2
=8(A^2)^2_{ij}. \]
Since $A$ is strictly upper triangular, $(A^2)_{ij}\neq0$ only if $i+1<j$, so
\[ \sum_{\gamma\in\Phi_2}S_\gamma\leq8\sum_{i+1<j}(A^2)^2_{ij}=8||A^2||^2\leq8||A||^4=8, \]
as we wished. 

\subsection{The $D_n$ family}
Recall that $\Delta^+=\Delta_1^+\cup\Delta_2^+$
where $\Delta_1^+=\{\theta_i-\theta_j|i<j\}$ and
$\Delta_2^+=\{\theta_i+\theta_j|i<j\}$. Consider the disjoint
union $\Phi=\Phi_1\cup\Phi_2\cup\Phi_3$ where
$\Phi_1=\Phi\cap2\Delta_1^+$, $\Phi_2=\Phi\cap(\Delta_1^++\Delta_2^+)$
and $\Phi_3=\Phi\cap2\Delta^+_2$.

Let $A=(a_{ij})$ be as in case $A_{n-1}$ and define in addition
the $n\times n$ matrix $B=(b_{ij})$ with
\[ b_{ij}=\left\{\begin{array}{cl}r_{\theta_i+\theta_j}&\mbox{if $i<j$,}\\
0 &\mbox{otherwise.}
\end{array}\right. \]

Similar to case $A_{n-1}$, we have
\[ \sum_{\gamma\in\Phi_1}S_\gamma\leq10||A||^4. \]
We shall show that
\begin{equation}\label{f2f3}
  \sum_{\gamma\in\Phi_2}S_\gamma\leq56||A||^2||B||^2\quad\mbox{and}
  \quad\sum_{\gamma\in\Phi_3}S_\gamma\leq18||B||^4.
  \end{equation}
Using $||A||^2+||B||^2=||fp||^2=1$, 
it then easily follows that
\begin{eqnarray}\label{22} \nonumber
  \sum_{\gamma\in\Phi}S_\gamma&\leq&10||A||^4+56||A||^2\,||B||^2+18||B||^4\\ \nonumber
  &=& 10+4(9||B||-7||B||^2)\\
  &<&22 
\end{eqnarray}
and hence, owing to~(\ref{II}), $||\sff_p||_\infty\leq2\sqrt6$. 

An element $\gamma\in\Phi_2$ has one of the following forms:
\begin{itemize}
\item[(i)] $\gamma=\theta_i+\theta_j-\theta_\ell+\theta_k$ where $i<j<\ell<k$; or
\item[(ii)] $\gamma=\theta_i+\theta_j+\theta_k-\theta_\ell$ 
where $i<j<k<\ell$; or
\item[(iii)] $\gamma=\theta_i+\theta_j$ where $i\leq j$.
\end{itemize}
Therefore
\begin{equation}\label{eq:dn}
  \sum_{\gamma\in\Phi_2}S_\gamma=\sum_{i<j<\ell<k}S_{\theta_i+\theta_j-\theta_\ell+\theta_k}
+\sum_{i<j<k<\ell}S_{\theta_i+\theta_j+\theta_k-\theta_\ell}
+\sum_{i\leq j}S_{\theta_i+\theta_j}.
\end{equation}

In case~(i), the possible decompositions
$\gamma=\alpha+\beta$ are given, up to a permutation, by
$\alpha=\theta_i-\theta_\ell$ and $\beta=\theta_j+\theta_k$ or
$\alpha=\theta_i+\theta_k$ and $\beta=\theta_j-\theta_\ell$.
Here $m_{\alpha,\beta}=1$ so
\[ S_{\theta_i+\theta_j-\theta_\ell+\theta_k}=8a_{i\ell}b_{jk}b_{ik}a_{j\ell}. \]

In case~(ii), the possible decompositions
$\gamma=\alpha+\beta$ are given, up to a permutation, by
$\alpha=\theta_i+\theta_j$ and $\beta=\theta_k-\theta_\ell$ or
$\alpha=\theta_i+\theta_k$ and $\beta=\theta_j-\theta_\ell$ or 
$\alpha=\theta_i-\theta_\ell$ and $\beta=\theta_j+\theta_k$.
Here also $m_{\alpha,\beta}=1$ so
\[ S_{\theta_i+\theta_j+\theta_k-\theta_\ell}=8(b_{ij}a_{k\ell}b_{ik}a_{j\ell}
+b_{ij}a_{k\ell}a_{i\ell}b_{jk}+b_{ik}a_{j\ell}a_{i\ell}b_{jk}). \]

It follows that the sum of the first two terms in the right
hand-side of~(\ref{eq:dn}) is equal to
\begin{eqnarray}\nonumber
  \lefteqn{8\left(\sum_{i<j<k<\ell}b_{ij}a_{k\ell}b_{ik}a_{j\ell}+b_{ij}a_{k\ell}a_{i\ell}b_{jk}+\sum_{i<j<k,\ell}a_{i\ell}b_{jk}b_{ik}a_{j\ell}\right)}\\ \nonumber
  &\leq&8(\langle B,BAA^t\rangle+
  \langle B^t,BAA^t\rangle+
  \langle B^t,B^tAA^t\rangle)\\ \label{dn1}
  &\leq&24 ||A||^2||B||^2.   \\    \nonumber
\end{eqnarray}

To deal with the remaining term in~(\ref{eq:dn}),
note that in case~(iii) we can have
$\gamma=\alpha+\beta=2\theta_i$, which can happen only if
$\alpha=\theta_i-\theta_k$ and $\beta=\theta_i+\theta_k$
with $i<k$, up to a permutation. It follows that $m_{\alpha,\beta}=1$ and
\begin{equation}\label{iii}
 S_{2\theta_i}\leq\left(\sum_k2a_{ik}b_{ik}\right)^2=4(AB^t)_{ii}^2. 
\end{equation}
We can also have $\gamma=\alpha+\beta=\theta_i+\theta_j$ with $i<j$,
which can happen only if
$\alpha=\theta_i-\theta_k$ and $\beta=\theta_j+\theta_k$
with $i<k\neq j$, or
$\alpha=\theta_i+\theta_k$ and $\beta=\theta_j-\theta_k$
with $j<k$, up to a permutation. Note that $m_{\alpha,\beta}=\sqrt2$
so (we use $b_{jk}+b_{kj}=r_{\theta_j+\theta_k}$ for $j\neq k$)
\[ S_{\theta_i+\theta_j}\leq 8\left(\sum_ka_{ik}(b_{kj}+b_{jk})
+\sum_kb_{ik}a_{jk}\right)^2=8(A(B+B^t)+BA^t)^2_{ij}. \]
Note that $M:=(A(B+B^t)+(B+B^t)A^t)^2$ is a symmetric matrix
which majorates both $(AB^t)^2$ and
$(A(B+B^t)+BA^t)^2$ term by term. It follows that
\begin{eqnarray}\nonumber
  \sum_{i\leq j}S_{\theta_i+\theta_j}&\leq&4\sum_iM_{ii}
  +8\sum_{i<j}M_{ij}\\ \nonumber
  &=&4\sum_iM_{ii}
  +4\sum_{i\neq j}M_{ij}\\ \nonumber
  &=&4||A(B+B^t)+(B+B^t)A^t||^2\\ \nonumber
  &\leq&16||A(B+B^t)||^2\\ \label{dn2}
  &\leq&32||A||^2||B||^2.\\ \nonumber
\end{eqnarray}
The first estimate in~(\ref{f2f3}) now follows from~(\ref{dn1})
and~(\ref{dn2}). 

There remains to obtain the second estimate in~(\ref{f2f3}).
Suppose $\gamma\in\Phi_3$.
Two elements in $\Delta_2^+$ cannot form an obtuse angle, so
$\gamma$ can be written as a sum of two positive roots 
in more than two ways. It follows that
$\gamma=\theta_i+\theta_j+\theta_k+\theta_\ell$ with
$i<j<k<\ell$. The allowed decompositions $\gamma=\alpha+\beta$
are given, up to a permutation, by
$\alpha=\theta_i+\theta_j$ and $\beta=\theta_k+\theta_\ell$ or
$\alpha=\theta_i+\theta_k$ and $\beta=\theta_j+\theta_\ell$ or
$\alpha=\theta_i+\theta_\ell$ and $\beta=\theta_j+\theta_k$.
Here $m_{\alpha,\beta}=1$ so
\begin{equation}\label{dn-f3}
 S_{\theta_i+\theta_j+\theta_k+\theta_\ell}=8(b_{ij}b_{k\ell}b_{ik}b_{j\ell}+
b_{ij}b_{k\ell}b_{i\ell}b_{jk}+b_{ik}b_{j\ell}b_{i\ell}b_{jk}). 
\end{equation}
The last term on the right hand-side of this equation
is symmetric in~$i$ and $j$ and in~$k$ and~$\ell$, hence
\begin{eqnarray}\nonumber
  \sum_{\gamma\in\Phi_3}S_\gamma&=&8\sum_{i<j<k<\ell}b_{ij}b_{k\ell}b_{ik}b_{j\ell}+
  b_{ij}b_{k\ell}b_{i\ell}b_{jk}+b_{ik}b_{j\ell}b_{i\ell}b_{jk}\\ \nonumber
  &\leq&8\sum_{i<j<k<\ell}b_{ij}b_{k\ell}b_{ik}b_{j\ell}+
  b_{ij}b_{k\ell}b_{i\ell}b_{jk}+2\sum_{i,j,k,\ell}
  b_{ik}b_{j\ell}b_{i\ell}b_{jk}\\ \nonumber
  &\leq&8\langle B,B^2B^t\rangle+8\langle B,B(B^t)^2\rangle
  +2\langle B,BB^tB\rangle\\ \label{dn-f3-end}
  &\leq&18||B||^4,
\end{eqnarray}
as desired. 

\subsection{The $C_n$ family}

We have $\Delta^+=\Delta^+_1\cup\Delta_2^+$
where $\Delta_1^+=\{\theta_i-\theta_j|i<j\}$ and
$\Delta_2^+=\{\theta_i+\theta_j|i\leq j\}$. Consider the disjoint
union $\Phi=\Phi_1\cup\Phi_2\cup\Phi_3$ where
$\Phi_1=\Phi\cap2\Delta_1^+$, $\Phi_2=\Phi\cap(\Delta_1^++\Delta_2^+)$
and $\Phi_3=\Phi\cap2\Delta^+_2$.

Let $A$ and $B$ be as in case~$D_n$ except that now we set
\[ b_{ii}=r_{2\theta_i} \]
for $i=1,\ldots,n$.
Note that $||A||^2+||B||^2=||fp||^2=1$.

Similar to case $A_{n-1}$, we have
\[ \sum_{\gamma\in\Phi_1}S_\gamma\leq10||A||^4. \]
We shall show that
\begin{equation}\label{f2f3bis}
  \sum_{\gamma\in\Phi_2}S_\gamma\leq56||A||^2||B||^2\quad\mbox{and}
  \quad\sum_{\gamma\in\Phi_3}S_\gamma\leq18||B||^4.
  \end{equation}
It will then follow as in the case $D_n$ that
$\sum_{\gamma\in\Phi}S_\gamma\leq22$
and hence, owing to~(\ref{II}), $||\sff_p||_\infty\leq2\sqrt6$. 

An element $\gamma\in\Phi_2$ has one of the following forms:
\begin{itemize}
\item[(i)] $\gamma=\theta_i+\theta_j-\theta_\ell+\theta_k$ where $i<j<\ell<k$; or
\item[(ii)] $\gamma=\theta_i+\theta_j+\theta_k-\theta_\ell$ 
where $i<j<k<\ell$; or
\item[(iii)] $\gamma=\theta_i+2\theta_j-\theta_\ell$ where $i<j<\ell$; or
\item[(iv)] $\gamma=2\theta_i+\theta_k-\theta_\ell$ where $i<k<\ell$; or
  \item[(v)] $\gamma=\theta_i+\theta_j$ where $i\leq j$.
\end{itemize}
Therefore 
\begin{eqnarray}\nonumber
  \lefteqn{\sum_{\gamma\in\Phi_2}S_\gamma=\sum_{i<j<\ell<k}S_{\theta_i+\theta_j-\theta_\ell+\theta_k}+\sum_{i<j<k<\ell}S_{\theta_i+\theta_j+\theta_k-\theta_\ell} }\\ \label{eq:cn}
&&\qquad+\sum_{i<j<\ell}S_{\theta_i+2\theta_j-\theta_\ell}
+\sum_{i<k<\ell}S_{2\theta_i+\theta_k-\theta_\ell}
+\sum_{i\leq j}S_{\theta_i+\theta_j}.\\ \nonumber
\end{eqnarray}

Similar to the case of $D_n$, the 
sum of the first two terms in the right
hand-side of ~(\ref{eq:cn}) is equal to
\begin{equation}\label{cn-1-2}
  8\left(\sum_{i<j<k<\ell}b_{ij}a_{k\ell}b_{ik}a_{j\ell}+\underbrace{b_{ij}a_{k\ell}a_{i\ell}b_{jk}}_{(\ast)}+\sum_{i<j<k,\ell}a_{i\ell}b_{jk}b_{ik}a_{j\ell}\right)
\end{equation}
We now add the third and fourth terms of~(\ref{eq:cn}) to the second
term (marked $(\ast)$) in~(\ref{cn-1-2}) as follows. 
In case~(iii), the possible decompositions
$\gamma=\alpha+\beta$ are given, up to permutation, by 
$\alpha=\theta_i+\theta_j$ and $\beta=\theta_j-\theta_\ell$ 
or $\alpha=\theta_i-\theta_\ell$ and $\beta=2\theta_j$. In both cases
$m_{\alpha,\beta}=1$ so 
\begin{equation}\label{expr-1}
 S_{\theta_i+2\theta_j-\theta_\ell}=8b_{ij}a_{j\ell}a_{i\ell}b_{jj}. 
\end{equation}
Likewise we see in regard to~(iv) that 
\begin{equation}\label{expr-2} 
S_{2\theta_i+\theta_k-\theta_\ell}=8b_{ii}a_{k\ell}b_{ik}a_{i\ell}. 
\end{equation}
Note that~(\ref{expr-1}) and~(\ref{expr-2}) coincide 
with~$(\ast)$  if $j=k$, resp., $j=i$. 
We deduce that the sum of the first four terms in the right
hand-side of~(\ref{eq:cn}) is bounded above by
\begin{equation}\label{cn1}
8\sum_{i,j,k,\ell}b_{ij}a_{k\ell}b_{ik}a_{j\ell}+
b_{ij}a_{k\ell}a_{i\ell}b_{jk}+
b_{ik}a_{j\ell}a_{i\ell}b_{jk}\leq24||A||^2||B||^2
\end{equation}
similar to~(\ref{dn1}).  

Next we address the remaining term in~(\ref{eq:cn})
coming from case~(v).  
A reasoning analogous to~(\ref{iii}) yields 
\begin{equation*}
 S_{2\theta_i}\leq\left(\sum_k2a_{ik}b_{ik}\right)^2=4(AB^t)_{ii}^2. 
\end{equation*}
If $\gamma=\alpha+\beta=\theta_i+\theta_j$ with $i<j$,
then
$\alpha=\theta_i-\theta_k$ and $\beta=\theta_j+\theta_k$
with $i<k\neq j$, or
$\alpha=\theta_i+\theta_k$ and $\beta=\theta_j-\theta_k$
with $j<k$, up to a permutation, or 
$\alpha=\theta_i-\theta_j$ and $\beta=2\theta_j$. 
It follows that 
\[ S_{\theta_i+\theta_j}\leq 8\left(\sum_ka_{ik}(b_{kj}+b_{jk})
+\sum_kb_{ik}a_{jk}+a_{ij}b_{jj}\right)^2. \]
Let $\tilde B$ the $n\times n$ matrix with coefficients 
$\tilde B_{ij}=b_{ij}+b_{ji}$ if $i\neq j$ and $\tilde B_{ii}=b_{ii}$. 
Then $S_{\theta_i+\theta_j}\leq8(A\tilde B+BA^t)^2_{ij}$ and, 
as in~(\ref{dn2}), we obtain
\begin{equation}\label{cn2}
\sum_{i\leq j}S_{\theta_i+\theta_j}\leq 16||A||^2||\tilde B||^2\leq32||A||^2||B||^2.
\end{equation}
The first estimate in~(\ref{f2f3bis}) now follows from~(\ref{cn1})
and~(\ref{cn2}). 

There remains to obtain the second estimate in~(\ref{f2f3bis}).
Suppose $\gamma\in\Phi_3$.
Two elements in $\Delta_2$ cannot form an obtuse angle, so
$\gamma$ can be written as a sum of two positive roots 
in more than two ways. It follows that $\gamma$ has one of the 
following forms:
\begin{itemize}
\item[(i)] $\gamma=\theta_i+\theta_j+\theta_k+\theta_\ell$ with
$i<j<k<\ell$; or
\item[(ii)] $\gamma=2\theta_i+\theta_k+\theta_\ell$ with
$i<k<\ell$; or
\item[(iii)] $\gamma=\theta_i+2\theta_j+\theta_\ell$ with
$i<j<\ell$; or
\item[(iv)] $\gamma=\theta_i+\theta_j+2\theta_k$ with
$i<j<k$; or
\item[(v)] $\gamma=2\theta_i+2\theta_k$ with
$i<k$. 
\end{itemize}
Regarding case~(i), similar to~(\ref{dn-f3}) we compute
\begin{eqnarray} \nonumber
\lefteqn{ \sum_{i<j<k<\ell}S_{\theta_i+\theta_j+\theta_k+\theta_\ell}=}\\ \label{cn-f3}
&&8\sum_{i<j<k<\ell}
b_{ij}b_{k\ell}b_{ik}b_{j\ell}+\underbrace{8\sum_{i<j<k<\ell}
b_{ij}b_{k\ell}b_{i\ell}b_{jk}}_{(\ast\ast)}+8\sum_{i<j<k<\ell}b_{ik}b_{j\ell}b_{i\ell}b_{jk}. 
\end{eqnarray}
We now add the sums corresponding to cases~(ii), (iii), (iv) and~(v)
to~$(\ast\ast)$. In case~(ii), the possible decompositions
$\gamma=\alpha+\beta$ are given, up to permutation, by 
$\alpha=2\theta_i$ and $\beta=\theta_k+\theta_\ell$ 
or $\alpha=\theta_i+\theta_k$ and $\beta=\theta_i+\theta_\ell$. In both cases
$m_{\alpha,\beta}=1$ so 
\begin{equation}\label{cn-ii}
 S_{2\theta_i+\theta_k+\theta_\ell}=8b_{ii}b_{k\ell}b_{ik}b_{i\ell}; 
\end{equation}
this is~$(\ast\ast)$ with $i=j$. 
Cases~(iii) and~(iv) are analogous and we obtain
\begin{equation}\label{cn-iii-iv} 
S_{\theta_i+2\theta_j+\theta_\ell}=8b_{ij}b_{j\ell}b_{i\ell}b_{jj}\quad\mbox{and}\quad
S_{\theta_i+\theta_j+2\theta_k}=8b_{ij}b_{kk}b_{ik}b_{jk}, 
\end{equation}
which are~$(\ast\ast)$ with $j=k$ and $k=\ell$, respectively. 
In case~(v), the possible decompositions
$\gamma=\alpha+\beta$ are given, up to permutation, by 
$\alpha=2\theta_i$ and $\beta=2\theta_k$ 
or $\alpha=\beta=\theta_i+\theta_k$. Since $m_{2\theta_i,2\theta_k}=1$ and 
$m_{\theta_i+\theta_k,\theta_i+\theta_k}=\sqrt2$, we get
\begin{equation}\label{cn-v} 
S_{2\theta_i+2\theta_k}=4\sqrt2b_{ii}b_{kk}b_{ik}^2\leq8b_{ii}b_{kk}b_{ik}^2, 
\end{equation}
namely, $(\ast\ast)$ with $i=j$ and $k=\ell$. 
Collecting the estimates~(\ref{cn-f3}), (\ref{cn-ii}),
(\ref{cn-iii-iv}) and (\ref{cn-v}), we arrive at
\[ \sum_{\gamma\in\Phi_3}S_\gamma\leq8\sum_{i<j<k<\ell}b_{ij}b_{k\ell}b_{ik}b_{j\ell}
+8\sum_{i,j,k,\ell}b_{ij}b_{k\ell}b_{i\ell}b_{jk}
+8\sum_{i<j<k<\ell} b_{ik}b_{j\ell}b_{i\ell}b_{jk}. \]
Repeating the arguments in~(\ref{dn-f3-end}), we finally deduce
the second estimate in~(\ref{f2f3bis}). 

\subsection{The $B_n$ family}
We have $\Delta^+=\Delta^+_1\cup\Delta_2^+\cup\Delta^+_3$
where $\Delta_1^+=\{\theta_i-\theta_j|i<j\}$,
$\Delta_2^+=\{\theta_i+\theta_j|i<j\}$
and $\Delta_3^+=\{\theta_i|1\leq i\leq n\}$.
Consider the disjoint
union $\Phi=\Phi_1\cup\Phi_2\cup\Phi_3\cup\Phi_4\cup\Phi_5$ where
$\Phi_1=\Phi\cap2\Delta_1^+$, $\Phi_2=\Phi\cap(\Delta_1^++\Delta_2^+)$,
$\Phi_3=\Phi\cap2\Delta^+_2$, $\Phi_4=\Phi\cap(\Delta^+_1+\Delta_3^+)$
and $\Phi_5=\Phi\cap(\Delta^+_2+\Delta_3^+)$ (note that 
$2\Delta_3^+\subset(\Delta_1^++\Delta_2^+)$). 

Let $A$ and $B$ be as in case $D_n$ and set $C$ to be the column-vector
with $c_i=r_{\theta_i}$ as coefficients. Note that 
$||A||^2+||B||^2+||C||^2=1$. 

Reasoning as in the case of $D_n$, we see that 
\[ \sum_{\gamma\in\Phi_1\cup\Phi_3}S_\gamma\leq10||A||^4+18||B||^4. \]
We shall show that 
\begin{eqnarray} \label{bn-f2} 
 \sum_{\gamma\in\Phi_2}S_\gamma&\leq&56||A||^2||B||^2+16||A||\,||B||\,||C||^2
+2||C||^4 \\ \label{bn-f4}
\sum_{\gamma\in\Phi_4}S_\gamma&\leq&12||A||^2||C||^2 \\ \label{bn-f5}
\sum_{\gamma\in\Phi_5}S_\gamma&\leq&16||B||^2||C||^2 \\ \nonumber
\end{eqnarray}
Writing $X=||A||$, $Y=||B||$, $Z=||C||$, we then have
\begin{equation}\label{pq}
 \sum_{\gamma\in\Phi}S_\gamma\leq p(X,Y) +Z^2q(X,Y) +2Z^4, 
\end{equation}
where 
\[ p(X,Y)=10X^4+56X^2Y^2+18Y^4\quad\mbox{and}\quad
q(X,Y)=12X^2+16XY+16Y^2 \]
subject to $X^2+Y^2+Z^2=1$. 
Since $p$ is homogeneous of degree~$4$ and $X^2+Y^2=1-Z^2$,
we deduce from~(\ref{22}) that $p(X,Y)\leq22(1-Z^2)^2$. Further,
\[ q(X,Y)\leq16(X^2+XY+Y^2)\leq16\frac32(X^2+Y^2)=24(1-Z^2). \]
Together with~(\ref{pq}), this says
\[ \sum_{\gamma\in\Phi}S_\gamma\leq 22-20Z^2\leq22 \]
and hence $||\sff_p||_\infty\leq2\sqrt6$. 

An element $\gamma\in\Phi_2$ has one of the following forms:
\begin{itemize}
\item[(i)] $\gamma=\theta_i+\theta_j-\theta_\ell+\theta_k$ where $i<j<\ell<k$; or
\item[(ii)] $\gamma=\theta_i+\theta_j+\theta_k-\theta_\ell$ 
where $i<j<k<\ell$; or
\item[(iii)] $\gamma=\theta_i+\theta_j$ where $i\leq j$.
\end{itemize}
Therefore
\begin{equation}\label{eq:bn-1}
  \sum_{\gamma\in\Phi_2}S_\gamma=\sum_{i<j<\ell<k}S_{\theta_i+\theta_j-\theta_\ell+\theta_k}
+\sum_{i<j<k<\ell}S_{\theta_i+\theta_j+\theta_k-\theta_\ell}
+\sum_{i\leq j}S_{\theta_i+\theta_j}.
\end{equation}
The situation for the first two terms on the right hand-side 
of~(\ref{eq:bn-1}) is completely analogous to the case~$D_n$ and thus
we already know that 
\begin{equation}\label{eq:bn-2}
\sum_{i<j<\ell<k}S_{\theta_i+\theta_j-\theta_\ell+\theta_k}
+\sum_{i<j<k<\ell}S_{\theta_i+\theta_j+\theta_k-\theta_\ell}\leq24||A||^2||B||^2. 
\end{equation}

In case~(iii), the possible decompositions
$\gamma=\alpha+\beta$ are the same as those listed in the corresponding 
case for $D_n$ plus $\alpha=\theta_i$ and $\beta=\theta_j$, up to a 
permutation. We deduce that 
\begin{equation*}
 S_{2\theta_i}\leq\left(\sqrt2c_i^2+\sum_k2a_{ik}b_{ik}\right)^2
=2(\sqrt2AB^t+CC^t)_{ii}^2 
\end{equation*}
and 
\begin{eqnarray*}
 S_{\theta_i+\theta_j}&\leq&\left(2c_ic_j+2\sqrt2\sum_ka_{ik}(b_{kj}+b_{jk})
+2\sqrt2\sum_kb_{ik}a_{jk}\right)^2\\ 
&=&4(\sqrt2A(B+B^t)+\sqrt2BA^t+CC^t)^2_{ij} 
\end{eqnarray*}
for~$i<j$. We proceed similarly to~(\ref{dn2}) to estimate
\begin{eqnarray}\nonumber
\sum_{i\leq j}S_{\theta_i+\theta_j}&\leq&2||(\sqrt2A(B+B^t)+\sqrt2(B+B^t)A^t+CC^t)||^2\\ \nonumber
&\leq&2(2\sqrt2||A(B+B^t)||+||C||^2)^2\\ \nonumber
&\leq&2(4||A||\,||B||+||C||^2)^2\\ \label{eq:bn-4}
&\leq&32||A||^2||B||^2+16||A||\,||B||\,||C||^2+2||C||^4.
\end{eqnarray} 
The estimate~(\ref{bn-f2}) now follows from~(\ref{eq:bn-2}) 
and~(\ref{eq:bn-4}). 

An element $\gamma\in\Phi_4$ has one of the following forms:
\begin{itemize}
\item[(i)] $\gamma=\theta_i+\theta_j-\theta_k$ where $i<j<k$; or
\item[(ii)] $\gamma=\theta_i$. 
\end{itemize}
In case~(i), the possible decompositions
$\gamma=\alpha+\beta$ are given, up to a permutation, by
$\alpha=\theta_i-\theta_k$ and $\beta=\theta_j$ or
$\alpha=\theta_j-\theta_k$ and $\beta=\theta_i$. Using symmetry in~$i$ and~$j$
for the first inequality, we get
\begin{eqnarray}\nonumber
\sum_{i<j<k}S_{\theta_i+\theta_j-\theta_k}&=&8\sum_{i<j<k}a_{ik}c_ja_{jk}c_i\\ \nonumber
&\leq&4\sum_{i,j,k}a_{ik}c_ja_{jk}c_i\\ \label{eq:bn-5}
&\leq&4||A||^2||C||^2.
\end{eqnarray}

In case~(ii), the possible decomposition
$\gamma=\alpha+\beta$ is given, up to a permutation, by
$\alpha=\theta_i-\theta_j$ and $\beta=\theta_j$ with $i<j$. 
Therefore
\begin{eqnarray}\nonumber
\sum_iS_{\theta_i}&\leq&\sum_i\left(2\sqrt2\sum_ja_{ij}c_j\right)^2\\ \label{eq:bn-6}
 &\leq&8||A||^2||C||^2.
\end{eqnarray}
The estimate~(\ref{bn-f4}) now follows from~(\ref{eq:bn-5}) 
and~(\ref{eq:bn-6}). 

An element $\gamma\in\Phi_5$ must be of the following form
$\gamma=\theta_i+\theta_j+\theta_k$ with $i<j<k$. 
The possible decompositions
$\gamma=\alpha+\beta$ are given, up to a permutation, by
$\alpha=\theta_i+\theta_j$ and $\beta=\theta_k$ or
$\alpha=\theta_i+\theta_k$ and $\beta=\theta_j$ or
$\alpha=\theta_j+\theta_k$ and $\beta=\theta_i$, so
\[ S_{\theta_i+\theta_j+\theta_k}=8(b_{ij}c_kb_{ik}c_j+b_{ij}c_kb_{jk}c_i
+b_{ik}c_jb_{jk}c_i). \]
Using that the first term on the righthand side of the above 
expression is symmetric in~$k$ and $j$ and the last term 
is symmetric in~$i$ and $j$, we can write
\begin{eqnarray*}
\sum_{\gamma\in\Phi_5}S_\gamma&=&\sum_{i<j<k}S_{\theta_i+\theta_j+\theta_k} \\
&\leq&8\left(\frac12\sum_{i,j,k}b_{ij}c_kb_{ik}c_j+\sum_{i,j,k}b_{ij}c_kb_{jk}c_i
+\frac12\sum_{i,j,k}b_{ik}c_jb_{jk}c_i\right)\\
&\leq&16||B||^2||C||^2,
\end{eqnarray*}
which is precisely~(\ref{bn-f5}). 

\section{The case of non-complex representations}\label{real}

In this section we discuss
the case in which $\rho$ does not admit an invariant complex structure. 
Note that in this case
$\rho$ is a representation of $K$ on a real vector 
space $W$ whose complexification is a complex irreducible representation $V$
and thus amenable to the results proved in the previous sections.

Let $\sigma$ be the conjugation of $V$ over $W$. 
Denote the $K$-invariant inner product on $W$ by 
$\langle\cdot,\cdot\rangle_{\mathbb R}$; it naturally extends to a 
$K$-invariant Hermitian product $\langle\cdot,\cdot\rangle$ on $V$.
The real part of $\langle\cdot,\cdot\rangle$ is an inner product 
on $V$, also denoted by $\langle\cdot,\cdot\rangle_{\mathbb R}$,
such that the decomposition
$V=W\oplus iW$ into the $\pm1$-eigenspaces of $\sigma$ is 
$\langle\cdot,\cdot\rangle_{\mathbb R}$-orthogonal.

If $\rho$ has the same orbits as the isotropy representation of a symmetric
space, then Lemma~\ref{isoparametric} says that we can
choose~$p$ in the unit sphere of $W$ such that $||\sff_p||\leq2$, and this
number is smaller than ($\sqrt2$ times) the numbers listed
in Table~1, so hereafter we may assume that is not the case for~$\rho$. 

Let $v_\lambda$ be a unit highest
weight vector of $V$. The component of $v_\lambda$ in $W$ is
$w=\frac12(v_\lambda+v_{-\lambda})$
where $v_{-\lambda}=\sigma(v_\lambda)$ is a unit lowest weight vector.
We choose $p=\sqrt2 w$ as a unit vector in $W$. For the sake
of clarity, below we denote the second 
fundamental forms of the $K$-orbits through~$v_\lambda$ in $V$ and 
through~$p$ in $W$ respectively by $\sff^V_{v_\lambda}$ and $\sff^W_p$. 
We next show that
\begin{equation}\label{real-2}
  ||\sff_p^W||_\infty\leq\sqrt2\,||\sff^V_{v_\lambda}||_\infty;
\end{equation}
note that 
we already know how to estimate the latter norm. 

Choose $x\in\Lk$ such that $||xp||=1$ 
and $||\sff^W_p(xp,xp)||=||\sff^W_p||_\infty$, and
put 
\[ \xi=\sff^W_p(xp,xp)/||\sff^W_p(xp,xp)||\in\nu_p(Kp)\subset W. \]
The normal space to $Kp$ at~$p$ in the unit sphere of $W$
is contained in the normal space to $Kv_\lambda$ at $v_\lambda$ in the 
unit sphere of $V$, so we compute 
\begin{eqnarray}\nonumber
||\sff^W_p||_\infty&=&\langle xxp,\xi\rangle_{\mathbb R} \\  \nonumber
&=&\sqrt2\langle xxw,\xi\rangle_{\mathbb R} \\ \nonumber
&=&\sqrt2\langle xxv_\lambda,\xi\rangle_{\mathbb R} \\ \nonumber
&\leq&\sqrt2||\sff_{v_\lambda}^V(xv_\lambda,xv_\lambda)||\\ \label{real-1}
&\leq&\sqrt2||\sff^V_{v_\lambda}||_\infty||xv_\lambda||^2.
\end{eqnarray}
Now $w=\frac1{\sqrt2}p$ so
\begin{eqnarray*}
\frac12&=&||xw||^2\\
&=&\frac14(||xv_\lambda||^2+||xv_{-\lambda}||^2+2\Re\langle xv_\lambda,
xv_{-\lambda}\rangle)\\
&=&\frac12(||xv_\lambda||^2+\Re\langle xv_\lambda,
xv_{-\lambda}\rangle) 
\end{eqnarray*}
equals $\frac12||xv_\lambda||^2$ if 
$\langle\Lk v_\lambda,\Lk v_{-\lambda}\rangle=\{0\}$;
in this case $||xv_\lambda||=1$ and
hence~(\ref{real-2}) follows from~(\ref{real-1}).
It remains to show that, under our assumptions,   
the case $\langle\Lk v_\lambda,\Lk v_{-\lambda}\rangle\neq\{0\}$
cannot happen. 

\begin{lem}\label{modest}
  If  $\langle\Lk v_\lambda,\Lk v_{-\lambda}\rangle\neq\{0\}$,
  then $\rho$ has the same orbits as
  the isotropy representation of a symmetric space. 
\end{lem}

\Pf The assumption implies that $v_{-\lambda}\in\mathcal U^2(\Lg)v_\lambda$,
where $\mathcal U^2(\Lg)$ is the second level in the natural filtration
of the enveloping algebra of $\Lg=\Lk\otimes_{\mathbb R}\C$, so Dadok's invariant 
$k(\lambda)=2~$~\cite[Propositions~4.1 and~4.4]{gt}. All such 
representations are orbit equivalent to isotropy representations of 
symmetric spaces (page~212 in~\cite{gt}). \EPf

\section{Conclusion}

The proof of the Theorem follows from results in the previous
sections. In fact, in view of~(\ref{focal-radius}), it suffices
to establish $||\sff_p||_\infty\leq C$ or $C\sqrt2$,
according to whether $\rho$ admits an invariant complex structure or not, 
for some point $p\in S^n$. 

Assume first $\rho$ admits an invariant complex structure.
Due to Proposition~\ref{simple}
(and the fact that the numbers in Table~1 are bigger than $\sqrt2$),
we may assume the group $K$ is simple.
The estimates for $K$ of exceptional type are discussed
in Subsection~\ref{exceptional}, whereas those for classical $K$
are done in Section~\ref{classical}.

The case $\rho$ does not admit an invariant complex structure
follows from~(\ref{real-2}) and the previous case.
This finishes the proof of the Theorem.

\section{Appendix: Computing $C_\Delta$}
The values of $C_\Delta$ for exceptional root systems
given in subsection~\ref{exceptional}
can be computed using the SageMath code below.
The function \texttt{m\_sq} takes as input two roots, $\alpha$ and $\beta$,
and returns the value of $m_{\alpha, \beta}^2$. The function \texttt{C\_delta}
takes as input a Cartan type written as a string and computes $C_\Delta$ for
that Cartan type by cycling though all positive root pairs.

{\footnotesize
\begin{lstlisting}[language=python]
def m_sq(alpha, beta):
    if alpha == beta:
        return 2
    if alpha.dot_product(beta)>=0:
        return 1
    else:
        return 2
        
def C_delta(T):
    pos_roots =  RootSystem(T).ambient_space().positive_roots()
    C = 0
    sums = {}
    for alpha in pos_roots:
        for beta in pos_roots:
            if alpha + beta in sums.keys():
                sums[alpha + beta] = sums[alpha + beta] + m_sq(alpha, beta)
            else:
                sums[alpha + beta] = m_sq(alpha, beta)
            if sums[alpha + beta] > C:
                C = sums[alpha + beta]
    return(C)
\end{lstlisting}
}
\bibliographystyle{amsalpha}
\bibliography{ref}

\end{document}